\documentclass{amsart}
\usepackage{amssymb, latexsym}

\def\bi{\begin{itemize}}
\def\bs{\begin{split}}
\def\es{\end{split}}
\def\ba{\begin{align}}
\def\bas{\begin{align*}}
\def\ea{\end{align}}
\def\eas{\end{align*}}

\def\C{{\mathbb C}}
\def\R{{{\mathbb R}}}
\def\Z{{{\mathbb Z}}}
\def\N{{{\mathbb N}}}

\def\emph#1{{\it #1}}
\def\textbf#1{{\bf #1}}

\newcommand{\lrs}{\langle \nabla \rangle^s}
\newcommand{\lrn}{\langle \nabla \rangle}

\newcommand{\lnr}{\langle\nabla\rangle}
\newcommand{\lqlr}{{L^q_t L^r_x}}
\newcommand{\ir}{{I \times \R^d}}

\newcommand{\rr}{{\R \times \R^d}}
\newcommand{\lli}{{\lqlr (\ir)}}
\newcommand{\llr}{{\lqlr (\rr)}}

\newcommand{\eps}{{\varepsilon}}

\theoremstyle{plain}
\newtheorem{theorem}{Theorem}
\newtheorem{definition}[theorem]{Definition}

\newtheorem{proposition}[theorem]{Proposition}
\newtheorem{lemma}[theorem]{Lemma}
\newtheorem{corollary}[theorem]{Corollary}

\numberwithin{equation}{section} \numberwithin{theorem}{section}

\begin{document}

\title[Blowup for the $L^2$-critical focusing NLS below $H^1$]
{On the blowup for the $L^2$-critical focusing nonlinear Schr\"odinger equation in higher dimensions below the energy class}
\author{Monica Visan}
\address{Institute for Advanced Study, Princeton NJ 08540}
\email{mvisan@ias.edu}
\author{Xiaoyi Zhang}
\address{Academy of Mathematics and System Sciences, Chinese Academy of Sciences}
\email{zh.xiaoyi@gmail.com}

\vspace{-0.3in}
\begin{abstract}
We consider the focusing mass-critical nonlinear Schr\"odinger equation and prove that blowup solutions to this equation with initial data
in $H^s(\R^d)$, $s > s_0(d)$ and $d\geq 3$, concentrate at least the mass of the ground state at the blowup time.  This extends
recent work by J.~Colliander, S.~Raynor, C.~Sulem, and J. D.~Wright, \cite{crsw}, T.~Hmidi and S.~Keraani, \cite{hk}, and N.~Tzirakis, \cite{tzirakis},
on the blowup of the two-dimensional and one-dimensional mass-critical focusing NLS below the energy space to all dimensions $d\ge 3$.
\end{abstract}

\maketitle

\section[And now]{Introduction}

We consider the initial value problem for the focusing $L_x^2$-critical nonlinear Schr\"odinger equation
\begin{equation}\label{equation}
\begin{cases}
i u_t +\Delta u = -|u|^{\frac 4{d}}u\\
u(0,x) = u_0(x)\in H^s(\R^d), \, s\geq 0,
\end{cases}
\end{equation}
where $u(t,x)$ is a complex-valued function in spacetime $\R\times\R^d$, $d\ge 3$.

It is well known (see, for example, \cite{cazenave:book}) that the Cauchy problem \eqref{equation} is locally wellposed in $H^s(\R^d)$
for $s\geq 0$.  Moreover, the unique solution obeys conservation of \emph{mass}:
$$
M(u(t)):=\int_{\R^d} |u(t,x)|^2\, dx = M(u_0).
$$
If $s\geq 1$, the \emph{energy} is also finite and conserved:
$$
E(u(t)):=\int_{\R^d} \bigl(\tfrac 12 |\nabla u(t,x)|^2 -\tfrac{d}{2(d+2)} |u(t,x)|^{2+\frac{4}{d}}\bigr)\, dx = E(u_0).
$$

This equation has a natural scaling.  More precisely, the map
\begin{equation}\label{scaling}
u(t,x)\mapsto u_\lambda(t,x):=\lambda^{\frac d2}u(\lambda^2 t,\lambda x)
\end{equation}
maps a solution to \eqref{equation} to another solution to \eqref{equation}.  The reason why this equation is called
$L_x^2$-critical (or mass-critical) is because the scaling \eqref{scaling} also leaves the mass invariant.

Equation \eqref{equation} is subcritical for $s>0$.  In this case, \eqref{equation} is wellposed in $H^s(\R^d)$
and the lifespan of the local solution depends only on the $H^s_x$-norm of the initial data (see \cite{cazenave:book}).
Denote by $T^*>0$ the maximal forward time of existence.  As a consequence of the local well-posedness theory, we have
the following blowup criterion:
$$
\mbox{either } \quad T^*=\infty \quad \mbox{ or } \quad T^*<\infty \mbox{ and } \lim_{t\to T^*} \|u(t)\|_{H_x^s}=\infty.
$$

The blowup behavior for solutions from $H_x^1$ initial data has received a lot of attention.  The results are closely related
to the ground state $Q$ which is the unique positive radial solution to the elliptic equation
$$
\Delta Q-Q+|Q|^{\frac 4d}Q=0.
$$
Using the sharp Gagliardo-Nirenberg inequality (see \cite{weinstein:ground state}),
\begin{equation}\label{gn}
\|u\|_{L_x^{\frac 4d+2}}^{\frac 4d+2}\le C_d\|u\|_{L_x^2}^{\frac 4d}\|\nabla u\|_{L_x^2}^2 \quad \text{with}\quad C_d:=\frac {d+2}{d\,\|Q\|^{\frac 4d}_{L_x^2}},
\end{equation}
it is not hard to see that the mass of the ground state is the minimal mass required for the solution to develop
a singularity.  Indeed, in the case $\|u_0\|_{L_x^2}<\|Q\|_{L_x^2}$, \eqref{gn} combined with the conservation of energy imply that
the solution to \eqref{equation} is global.  This is sharp since the pseudoconformal invariance of the equation \eqref{equation}
allows us to built a solution with mass equal to that of the ground state that blows up at time $T^*$:
$$
u(t,x):=|T^*-t|^{-\frac d2}e^{[i(T^*-t)^{-1}-i|x|^2(T^*-t)^{-1}]}Q\Bigl(\frac x{T^*-t}\Bigr).
$$
Moreover, F.~Merle, \cite{merle1}, showed that up to the symmetries of \eqref{equation}, this is the only blowup solution with minimal mass.
Furthermore, any blowup solution must concentrate at least the mass of the ground state at the blowup time; more precisely, as shown in
\cite{merle-tsutsumi}, there exists $x(t)\in \R^d$ such that
\begin{align}\label{H^1 concentration}
\forall R>0,\ \lim_{t\to T^*} \int_{|x-x(t)|\le R}|u(t,x)|^2dx\ge \int_{\R^d} Q^2dx.
\end{align}

Of course, the goal is to establish all these properties for blowup solutions from data in $L_x^2$ rather than $H^1_x$.  Unfortunately,
all the methods used in the $H^1_x$ setting break down at the $L_x^2$ level.  Moreover, as \eqref{equation} is $L_x^2$-critical, even
the local well-posedness theory in $L_x^2$ is substantially different from that in $H^s_x$ for $s>0$.  Specifically, the lifespan
of the local solution depends on the profile of the initial data, rather than on its $L_x^2$-norm (see \cite{cazenave:book}).
In particular, this leads to the following blowup criterion:
$$
\mbox{either } \quad T^*=\infty \quad \mbox{ or } \quad T^*<\infty \mbox{ and } \|u\|_{L_{t,x}^{2+\frac 4d}([0,T^*)\times \R^d)}=\infty.
$$
From the global theory for small data (see \cite{cazenave:book}), we know that if the mass of the initial data is sufficiently small, then
there exists a unique global solution to \eqref{equation}.  However, for large (but finite) mass initial data, that is also sufficiently
smooth and decaying, the viriel identity guarantees that finite time blowup occurs; see, \cite{glassey, zak}.

The first blowup result for general $L_x^2$ initial data belongs to J.~Bourgain, \cite{bourgain:blowup}, who proved
the following parabolic concentration of mass at the blowup time:
\begin{align}\label{L^2 concentration}
\lim_{t\nearrow T^*} \sup_{\substack{\text{cubes}\ I\subset \R^2 \\ \text{side}(I)<(T^*-t)^{\frac 12}}}
\Bigl(\int_I |u(t,x)|^2\, dx\Bigr)^{\frac 12}\ge c(\|u_0\|_{L_x^2})>0,
\end{align}
where $c(\|u_0\|_{L_x^2})$ is a small constant depending on the mass of the initial data.  This result was extended to dimension $d=1$ by
S.~Keraani, \cite{keraani}, and to dimensions $d\geq 3$ by P.~Begout and A.~Vargas, \cite{begout:vargas}.  The conjecture is that rather that
$c(\|u_0\|_{L_x^2})$ on the right-hand side of \eqref{L^2 concentration}, one should have the mass of the ground state as in
\eqref{H^1 concentration}; however, this appears to be a very difficult problem.

The goal of this paper is to reproduce as much of the $H^1_x$ theory as we can at lower regularity (but unfortunately, well above $L^2_x$)
in dimensions three and higher.  In two dimensions, this was pursued by J.~Colliander, S.~Raynor, C.~Sulem, and J. D.~Wright, \cite{crsw}; using the \emph{I-method}
(introduced by J.~Colliander, M.~Keel, G.~Staffilani, H.~Takaoka, and T.~Tao, cf. \cite{ckstt:low1}), they established a mass concentration property
(at the blowup time) for radial blowup solutions from initial data in $H^s(\R^2)$ with $s>\frac{1+\sqrt{11}}5$:
\begin{align}\label{H^s concentration 2d}
\limsup_{t\nearrow T^*} \|u\|_{L^2(|x|<(T^*-t)^{\frac s2}\gamma(T^*-t))}\geq \|Q\|_{L_x^2},
\end{align}
where $\gamma(s)\to \infty$ arbitrarily slowly as $s\to 0$.  While the right-hand side in \eqref{H^s concentration 2d} is the conjectured one,
we should remark that the width of concentration, that is, $(T^*-t)^{\frac s2 +}$ is larger than the expected concentration width,
$(T^*-t)^{\frac 12 +}$.  Moreover, the $\limsup$ in \eqref{H^s concentration 2d} is not present in the $H^1_x$ theory (see \eqref{H^1 concentration});
the reason for its appearance is basically the lack of information on the blowup rate of the $H^s_x$-norm.

T.~Himidi and S.~Keraani, \cite{hk}, removed the radial assumption used in \cite{crsw}.  Moreover, they showed that, up to symmetries of the equation,
the ground state $Q$ is the profile for blowup solutions with minimal mass and initial data in $H^s(\R^2)$ with $s>\frac{1+\sqrt{11}}5$.
We will make use of the key innovation of their work, namely Lemma~\ref{compact} below.

In one dimension, recent work by N. Tzirakis, \cite{tzirakis}, established the analogue of \eqref{H^s concentration 2d} for arbitrary initial data
in $H^s(\R)$, $s>\frac{10}{11}$.

Our contribution is to treat (non-radial) data in $H^s(\R^d)$ for $s>s_0(d)$ where
\begin{equation}\label{s0d}
s_0(d):=\left\{
\begin{array}{cccc}
&\frac {1+\sqrt{13}}5, & \text{for}{\vrule width 0mm depth 2ex height 0mm} \quad d=3,\\
&\frac {8-d+\sqrt{9d^2+64d+64}}{2(d+10)}, & \text{for} \quad d\ge4.
\end{array}
\right.
\end{equation}
We prove

\begin{theorem}\label{blowup}
Assume $d\ge 3$ and $s>s_0(d)$.  Let $u_0\in H^s(\R^d)$ such that the corresponding solution $u$ to \eqref{equation} blows up at time $0<T^*<\infty$.
Then, there exists a function $V\in H_x^1$ such that $\|V\|_2\ge \|Q\|_2$ and there exist sequences
$\{t_n,\rho_n,x_n\}_{n\geq 1}\subset \R_+\times\R_+^*\times \R^d$ satisfying
$$
t_n\nearrow T^* \text{ as } n\to\infty \quad \text{and} \quad \rho_n\lesssim(T^*-t_n)^{\frac s2},\ \ \forall\, n\geq 1
$$
such that
$$
\rho_n^{\frac d2} u(t_n,\rho_n\cdot +x_n)\rightharpoonup  V \quad \text{weakly} \quad \text{as } n\to\infty.
$$
\end{theorem}


As a consequence of Theorem~\ref{blowup}, we establish the following mass concentration property for blowup solutions:

\begin{theorem}\label{concentration}
Assume $d\ge 3$ and $s>s_0(d)$. Let $u_0\in H^s(\R^d)$ such that the corresponding solution $u$ to \eqref{equation} blows up at time $0<T^*<\infty$.
Let $\alpha(t)>0$ be such that
$$
\lim_{t\nearrow T^*}\frac {(T^*-t)^{\frac s2}}{\alpha(t)}=0.
$$
Then, there exists $x(t)\in \R^d$ such that
$$
\limsup_{t\nearrow T^*}\int_{|x-x(t)|\le \alpha(t)}|u(x,t)|^2dx\ge \int_{\R^d}Q^2dx.
$$
\end{theorem}

Under the additional hypothesis that the mass of the initial data equals the mass of the ground state, we may upgrade Theorem~\ref{blowup}
to the following:

\begin{theorem}\label{gs blowup}
Assume $d\ge 3$ and $s>s_0(d)$.  Let $u_0\in H^s(\R^d)$ with $\|u_0\|_2=\|Q\|_2$ such that the corresponding solution $u$ to \eqref{equation}
blows up at time $0<T^*<\infty$.  Then, there exist sequences $\{t_n,\theta_n,\rho_n,x_n\}_{n\geq 1} \subset \R_+\times S^1\times \R_+^*\times \R^d$
satisfying
$$
t_n\nearrow T^* \text{ as } n\to\infty \quad \text{and} \quad \rho_n\lesssim(T^*-t_n)^{\frac s2},\ \ \forall\, n\geq 1
$$
such that
$$
\rho_n^{\frac d2} e^{i\theta_n} u(t_n,\rho_nx+x_n) \rightarrow Q \quad \text{strongly in } H_x^{\tilde s -},
$$
where
$$
\tilde s:=\frac {2d+8s+s^2d(2-\min\{1,\frac 4d\})}{4d+16s-s^2(d\min\{1,\frac 4d\}+8)}.
$$
\end{theorem}

In a nutshell, Theorem~\ref{gs blowup} says that up to the symmetries for \eqref{equation}, the ground state is the profile for blowup
solutions with minimal mass and initial data in $H^s_x$, $s>s_0(d)$.  Alas, we only show this is true along a sequence of times.

To prove Theorems \ref{blowup} through \ref{gs blowup}, we will rely on the \emph{I-method} and Lemma~\ref{compact}.
The idea behind the \emph{I-method} is to smooth out the initial data in order to access the theory available at $H^1_x$ regularity.
To this end, one introduces the Fourier multiplier $I$, which is the identity on low frequencies and behaves like a fractional integral operator
of order $1-s$ on high frequencies.  Thus, the operator $I$ maps $H_x^s$ to $H_x^1$.  However, even though we do have energy conservation for
\eqref{equation}, $Iu$ is not a solution to \eqref{equation} and hence, we expect an energy increment.  The key is to prove that on intervals of
local well-posedness, the modified energy $E(Iu)$ is an `almost conserved' quantity and grows much slower than the modified kinetic energy
$\|\nabla Iu\|_{L_x^2}^2$. This requires delicate estimates on the commutator between $I$ and the nonlinearity.  In dimensions one and two,
the nonlinearity is algebraic and one can write the commutator explicitly using the Fourier transform and control it by multilinear analysis
and bilinear estimates (see \cite{crsw, tzirakis}).  However, in dimensions $d\geq 3$ this method fails.  Instead, we will have to rely on more
rudimentary tools such as Strichartz and fractional chain rule estimates in order to control the commutator.

The remainder of this paper is organized as follows:  In Section~2, we introduce notation and prove some lemmas that will be useful.
In Section~3, we revisit the $H_x^s$ local well-posedness theory for \eqref{equation}.  Section~4 is devoted to controlling the modified
energy increment.  In Sections 5 through 7 we prove Theorems \ref{blowup} through \ref{gs blowup}.

After this work was submitted, we were informed of an independent paper attacking the same problem, \cite{sijia}.  However, there appear to be
several gaps in their argument, for example, in estimating (4.19).  While it seems possible to remedy their errors, this would result in a larger
value for $s_0(d)$ than that claimed in their paper, which is already inferior to that given here.  We contend that all overlap between
this paper and \cite{sijia} can be attributed to the precursors, \cite{crsw} and \cite{hk}.

\subsection*{Acknowledgement} The authors would like to thank the organizers of the Nonlinear Dispersive Equations program at M.S.R.I. during which
this work was completed.  The authors would also like to thank Jim Colliander for introducing us to the problem.  The second author was supported by the
NSF grant No. 10601060 (China).

%
%
%
%

\section{Preliminaries}
We will often use the notation $X \lesssim Y$ whenever there exists some constant $C$ so that $X \leq CY$. Similarly, we will use $X
\sim Y$ if $X \lesssim Y \lesssim X$.  We use $X \ll Y$ if $X \leq cY$ for some small constant $c$. The derivative operator $\nabla$
refers to the space variable only.  We use $A\pm$ to denote $A\pm\eps$ for any sufficiently small $\eps>0$; the implicit constant in an inequality
involving this notation is permitted to depend on $\eps$.

We use $L_x^r(\R^d)$ to denote the Banach space of functions
$f:\R^d\to \C$ whose norm
$$
\|f\|_r:=\Bigl(\int_{\R^d} |f(x)|^r dx\Bigr)^{\frac{1}{r}}
$$
is finite, with the usual modifications when $r=\infty$.

We use $\lqlr$ to denote the spacetime norm
$$
\|u\|_{q,r}:=\|u\|_{\llr} :=\Bigl(\int_{\R}\Bigl(\int_{\R^d}
|u(t,x)|^r dx \Bigr)^{q/r} dt\Bigr)^{1/q},
$$
with the usual modifications when either $q$ or $r$ are infinity, or
when the domain $\R \times \R^d$ is replaced by some smaller
spacetime region.  When $q=r$ we abbreviate $\lqlr$ by $L^q_{t,x}$.

We define the Fourier transform on $\R^d$ to be
$$
\hat f(\xi) := \int_{\R^d} e^{-2 \pi i x \cdot \xi} f(x) dx.
$$

We will make use of the fractional differentiation operators $|\nabla|^s$ defined by
$$
\widehat{|\nabla|^sf}(\xi) := |\xi|^s \hat f (\xi).
$$
These define the homogeneous Sobolev norms
$$
\|f\|_{\dot H^s_x} := \| |\nabla|^s f \|_{L^2_x}
$$
and the more common inhomogeneous Sobolev norms
$$
\|f\|_{H^{s,p}_x}:=\|\langle \nabla \rangle^sf\|_{L_x^p}, \quad \text{where} \quad \langle \nabla \rangle:=(1+|\nabla|^2)^{\frac 12}.
$$
We will often denote $H^{s,2}_x$ by $H^s_x$.

Let $F(z):=-|z|^{\frac{4}{d}}z$ be the function that defines the nonlinearity in \eqref{equation}.  Then,
$$
F_z(z):=\frac{\partial F}{\partial z}(z)=-\tfrac{2+d}{d}|z|^{\frac 4{d}} \quad \text{and} \quad
F_{\bar z}(z):=\frac{\partial F}{\partial \bar z}(z)=-\tfrac{2}d |z|^{\frac{4}{d}}\frac{z}{\bar z}.
$$
We write $F'$ for the vector $(F_z,F_{\bar z})$ and adopt the notation
$$
w\cdot F'(z):=wF_{z}(z)+\bar wF_{\bar z}(z).
$$
In particular, we observe the chain rule
$$
\nabla F(u)=\nabla u\cdot F'(u).
$$
Clearly $F'(z)=O(|z|^{\frac{4}{d}})$ and we have the H\"older continuity estimate
\begin{align}\label{holder continuity}
|F'(z)-F'(w)|\lesssim |z-w|^{\min\{1,\frac{4}{d}\}}(|z|+|w|)^{\frac{4}{d}-\min\{1,\frac{4}{d}\}}
\end{align}
for all $z, w\in \C$.  By the Fundamental Theorem of Calculus,
$$
F(z+w)-F(z)=\int_0^1 w\cdot F'(z+\theta w)d\theta
$$
and hence
$$
F(z+w)=F(z)+O(|w||z|^\frac{4}{d})+O(|w|^{\frac{4+d}{d}})
$$
for all complex values $z$ and $w$.

Let $e^{it\Delta}$ be the free Schr\"odinger propagator.  In physical space this is given by the formula
$$
e^{it\Delta}f(x) = \frac{1}{(4 \pi i t)^{d/2}} \int_{\R^d}
e^{i|x-y|^2/4t} f(y) dy
$$
for $t\neq 0$ (using a suitable branch cut to define $(4\pi
it)^{d/2}$), while in frequency space one can write this as
\begin{equation}\label{fourier rep}
\widehat{e^{it\Delta}f}(\xi) = e^{-4 \pi^2 i t |\xi|^2}\hat f(\xi).
\end{equation}
In particular, the propagator obeys the \emph{dispersive inequality}
\begin{equation}\label{dispersive ineq}
\|e^{it\Delta}f\|_{L^\infty_x} \lesssim
|t|^{-\frac{d}{2}}\|f\|_{L^1_x}
\end{equation}
for all times $t\neq 0$.

We also recall \emph{Duhamel's formula}
\begin{align}\label{duhamel}
u(t) = e^{i(t-t_0)\Delta}u(t_0) - i \int_{t_0}^t e^{i(t-s)\Delta}(iu_t + \Delta u)(s) ds.
\end{align}

\begin{definition}
A pair of exponents $(q,r)$ is called Schr\"odinger-\emph{admissible} if
$$
\frac{2}{q} +\frac{d}{r} = \frac{d}{2},\quad 2 \leq q,r \leq \infty, \quad \text{and} \quad (q,r,d)\neq (2,\infty,2).
$$
\end{definition}
Throughout this paper we will use the following admissible pairs:
$$
(2,\tfrac{2d}{d-2}) \quad \text{and} \quad (\gamma,\rho):=(\tfrac{2(d+2)}{d-2s},\tfrac{2d(d+2)}{d^2+4s}) \quad \text{with } 0<s<1.
$$
Let $\rho^*:=\frac{2(d+2)}{d-2s}$.  Using H\"older and Sobolev embedding, on any spacetime slab $\ir$ we estimate
\begin{align}\label{estimate gamma-rho}
\|F(u)\|_{\gamma',\rho'}
\lesssim |I|^{\frac{2s}{d}}\|u\|_{\gamma, \rho} \|u\|_{\gamma, \rho^*}^{\frac 4d}
\lesssim |I|^{\frac{2s}{d}}\|\langle \nabla\rangle ^s u\|_{\gamma, \rho}^{1+\frac 4d}.
\end{align}

Let $I \times \R^d$ be a spacetime slab; we define the Strichartz
norm
\begin{equation*}
\|u\|_{S^0(I)} := \sup_{(q,r)\, \text{admissible}} \|u \|_{\lli}.
\end{equation*}
We define the Strichartz space $S^0(I)$ to be the closure of all test functions
under the Strichartz norm $\|\cdot\|_{S^0(I)}$.  We use $N^0(I)$ to denote the dual space of $S^0(I)$.

We record the standard Strichartz estimates which we will invoke repeatedly throughout this paper (see \cite{strichartz, gv:Strichartz}
for $(q,r)$ admissible with $q>2$ and \cite{tao:keel} for the endpoint $(2,\tfrac{2d}{d-2})$):

\begin{lemma}\label{lemma linear strichartz}
Let $I$ be a compact time interval, $t_0\in I$, $s\geq 0$, and let $u$ be a solution to the forced Schr\"odinger equation
\begin{equation*}
i u_t + \Delta u =\sum_{i=1}^m F_i
\end{equation*}
for some functions $F_1$, $\cdots$, $F_m$.  Then,
\begin{equation}
\||\nabla| ^s u\|_{S^0(I)} \lesssim \|u(t_0)\|_{\dot H^s(\R^d)} +
\sum_{i=1}^m \||\nabla| ^s F_i\|_{L_t^{q_i'}L_x^{r_i'}
(I\times\R^d)}
\end{equation}
for any admissible pairs $(q_i,r_i)$, $1\leq i\leq m$.
\end{lemma}

We will also need some Littlewood-Paley theory.  Specifically, let $\varphi(\xi)$ be a smooth bump supported in the ball
$|\xi| \leq 2$ and equalling one on the ball $|\xi| \leq 1$.  For each dyadic number $N \in 2^\Z$ we define the
Littlewood-Paley operators
\begin{align*}
\widehat{P_{\leq N}f}(\xi) &:=  \varphi(\xi/N)\hat f (\xi),\\
\widehat{P_{> N}f}(\xi) &:=  [1-\varphi(\xi/N)]\hat f (\xi),\\
\widehat{P_N f}(\xi) &:=  [\varphi(\xi/N) - \varphi (2 \xi /N)] \hat f (\xi).
\end{align*}
Similarly we can define $P_{<N}$, $P_{\geq N}$, and $P_{M < \cdot \leq N} := P_{\leq N} - P_{\leq M}$, whenever $M$ and
$N$ are dyadic numbers.  We will frequently write $f_{\leq N}$ for $P_{\leq N} f$ and similarly for the other operators.
We recall the following standard Bernstein and Sobolev type inequalities:

\begin{lemma}\label{bernstein}
For any $1\le p\le q\le\infty$ and $s>0$, we have
\begin{align*}
\|P_{\geq N} f\|_{L^p_x} &\lesssim N^{-s} \| |\nabla|^s P_{\geq N} f \|_{L^p_x}\\
\| |\nabla|^s  P_{\leq N} f\|_{L^p_x} &\lesssim N^{s} \| P_{\leq N} f\|_{L^p_x}\\
\| |\nabla|^{\pm s} P_N f\|_{L^p_x} &\sim N^{\pm s} \| P_N f \|_{L^p_x}\\
\|P_{\leq N} f\|_{L^q_x} &\lesssim N^{\frac{d}{p}-\frac{d}{q}} \|P_{\leq N} f\|_{L^p_x}\\
\|P_N f\|_{L^q_x} &\lesssim N^{\frac{d}{p}-\frac{d}{q}} \| P_N
f\|_{L^p_x}.
\end{align*}
\end{lemma}

\vspace{0.4cm}

For $N>1$, we define the Fourier multiplier $I_N$ (cf. \cite{ckstt:low1}) by
$$
\widehat{I_N u}(\xi):=m_N(\xi)\hat u(\xi),
$$
where $m_N$ is a smooth radial decreasing function such that
$$
m_N(\xi)=\left\{
\begin{array}{cc}
1, \quad \text{if} \quad  |\xi|\le N\\
\bigl(\frac{|\xi|}{N}\bigr)^{s-1}, \quad \text{if} \quad |\xi|\ge 2N.
\end{array}
\right.
$$
Thus, $I_N$ is the identity operator on frequencies $|\xi|\le N$ and behaves like a fractional integral
operator of order 1-s on higher frequencies. In particular, $I_N$ maps $H^s_x$ to $H_x^1$; this allows us to
access the theory available for $H_x^1$ data.  We collect the basic properties of $I_N$ into the following:

\begin{lemma}\label{basic property}
Let $1<p<\infty$ and $0\leq\sigma\le s<1$.  Then,
\begin{align}
\|I_N f\|_p&\lesssim \|f\|_p \label{i1}\\
\||\nabla|^\sigma P_{> N}f\|_p&\lesssim N^{\sigma-1}\|\nabla I_N f\|_p \label{i2}\\
\|f\|_{H^s_x}\lesssim \|I_N f\|_{H^1_x}&\lesssim N^{1-s}\|f\|_{H^s_x}.\label{i3}
\end{align}
\end{lemma}
\begin{proof}
The estimate \eqref{i1} is a direct consequence of the multiplier theorem.

To prove \eqref{i2}, we write
$$
\||\nabla|^\sigma P_{> N} f\|_p=\|P_{> N}|\nabla |^\sigma(\nabla I_N)^{-1}\nabla I_N f\|_p.
$$
The claim follows again from the multiplier theorem.

Now we turn to \eqref{i3}.  By the definition of the operator $I_N$ and \eqref{i2},
\begin{align*}
\|f\|_{H^s_x}&\lesssim \|P_{\le N} f\|_{H^s_x}+\|P_{>N}f\|_2+\||\nabla|^s P_{>N} f\|_2\\
&\lesssim \|P_{\le N} I_N f\|_{H_x^1}+N^{-1}\|\nabla I_N f\|_2+N^{s-1}\|\nabla I_N f\|_2\\
&\lesssim \|I_Nf\|_{H^1_x}.
\end{align*}
On the other hand, since the operator $I_N$ commutes with $\langle \nabla\rangle^s$, we get
\begin{align*}
\|I_Nf\|_{H_x^1}
=\|\langle\nabla\rangle^{1-s}I_N\langle\nabla\rangle^s f\|_2
\lesssim N^{1-s}\|\langle\nabla\rangle^sf\|_2
\lesssim N^{1-s}\|f\|_{H^s_x},
\end{align*}
which proves the last inequality in \eqref{i3}.  Note that a similar argument also yields
\begin{align}\label{i4}
\|I_Nf\|_{\dot H^1_x}&\lesssim N^{1-s}\|f\|_{\dot H^s_x}.
\end{align}
\end{proof}

The estimate \eqref{i2} shows that we can control the high frequencies of a function $f$ in the Sobolev space $H_x^{\sigma,p}$ by the smoother function
$I_Nf$ in a space with a loss of derivative but a gain of negative power of $N$. This fact is crucial in extracting the negative power of $N$
when estimating the increment of the modified Hamiltonian.

In dimensions one and two, one can use multilinear analysis to understand commutator expressions like $F(I_Nu) - I_NF(u)$; on the Fourier side,
one can expand this commutator into a product of Fourier transforms of $u$ and $I_Nu$ and carefully measure the frequency interactions to derive
an estimate (see for example \cite{crsw, tzirakis}).  However, this is not possible in dimensions $d\geq 3$.  Instead, we will have to rely
on the following rougher (weaker, but more robust) lemma:

\begin{lemma}\label{bilinear}
Let $1<r,r_1,r_2<\infty$ be such that $\frac 1 r=\frac 1{r_1}+\frac 1{r_2}$ and let $0<\nu<s$.  Then,
\begin{align}\label{bilinear eq}
\|I_N(fg)-(I_Nf)g\|_r\lesssim N^{-(1-s+\nu)}\|I_Nf\|_{r_1}\|\langle\nabla\rangle^{1-s+\nu}g\|_{r_2}.
\end{align}
\end{lemma}
\begin{proof} Applying a Littlewood-Paley decomposition to $f$ and $g$, we write
\begin{align}
I_N(fg)-&(I_Nf)g \nonumber\\
&=I_N(fg_{\le 1})-(I_Nf)g_{\le1} +\sum_{1<M\in 2^{\Z}}\bigl[I_N(f_{\lesssim M}g_M)-(I_Nf_{\lesssim M})g_M \bigr] \nonumber\\
&\qquad \qquad +\sum_{1<M\in 2^{\Z}}\bigl[I_N(f_{\gg M}g_M)-(I_Nf_{\gg M})g_M\bigr] \nonumber\\
&=I_N(f_{\gtrsim N}g_{\le 1})-(I_Nf_{\gtrsim N}) g_{\le 1} +\sum_{N\lesssim M\in 2^{\Z}}\bigl[I_N(f_{\lesssim M} g_M) -(I_Nf_{\lesssim M})g_M\bigr] \nonumber\\
&\qquad\qquad +\sum_{1<M\in 2^{\Z}}\bigl[I_N(f_{\gg M}g_M)-(I_Nf_{\gg M})g_M \bigr]\nonumber\\
&=I +II +III.\label{difference}
\end{align}
The second equality above follows from the fact that the operator $I_N$ is the identity operator on frequencies $|\xi|\leq N$; thus,
\begin{align*}
I_N(f_{\ll N}g_{\le1})=(I_Nf_{\ll N})g_{\le 1} \quad \text{and} \quad I_N(f_{\lesssim M}g_{M})=(I_Nf_{\lesssim M})g_{M}& \quad \text{for all } M\ll N.
\end{align*}

We first consider $II$.  Dropping the operator $I_N$, by H\"older and Bernstein we estimate
\begin{align*}
\|I_N(f_{\lesssim M} g_M) -(I_Nf_{\lesssim M})g_M\|_r&\lesssim \|f_{\lesssim M}\|_{r_1}\|g_M\|_{r_2}\\
&\lesssim \bigl(\frac MN\bigr)^{1-s} \|I_Nf\|_{r_1}\|g_M\|_{r_2}\\
&\lesssim M^{-\nu}N^{-(1-s)}\|I_Nf\|_{r_1}\||\nabla|^{1-s+\nu}g\|_{r_2}.
\end{align*}
Summing over all $M$ such that $N\lesssim M\in 2^{\Z}$, we get
\begin{equation}\label{d1}
II\lesssim N^{-(1-s+\nu)}\|I_Nf\|_{r_1}\||\nabla|^{1-s+\nu}g\|_{r_2}.
\end{equation}

We turn now towards $III$.  Applying a Littlewood-Paley decomposition to $f$, we write each term in $III$ as
\begin{align*}
I_N(f_{\gg M}g_M)-(I_Nf_{\gg M})g_M
&=\sum_{1\ll k\in\N} \bigl[ I_N(f_{2^k M}g_M)-(I_Nf_{2^kM})g_M\bigr]\\
&=\sum_{\substack{1\ll k\in\N \\ N\lesssim 2^kM}} \bigl[I_N(f_{2^k M} g_M)-(I_Nf_{2^kM})g_M\bigr].
\end{align*}
To derive the second inequality, we used again the fact that the operator $I_N$ is the identity on frequencies $|\xi|\le N$.

We write
$$
[I_N(f_{2^k M} g_M)-(I_Nf_{2^kM})g_M]\widehat{\ }(\xi)=\int_{\xi=\xi_1+\xi_2}(m_N(\xi_1+\xi_2)-m_N(\xi_1))\widehat{f_{2^k M}}(\xi_1)\widehat{g_M}(\xi_2).
$$
For $|\xi_1|\sim 2^kM$, $k\gg 1$, and $|\xi_2|\sim M$, the Fundamental Theorem of Calculus implies
$$
|m_N(\xi_1+\xi_2)-m_N(\xi_1)|\lesssim 2^{-k} \bigl(\frac{2^kM}{N}\bigr)^{s-1}.
$$
By the Coifman-Meyer multilinear multiplier theorem, \cite{coifmey:1, coifmey:2}, and Bernstein, we get
\begin{align*}
\|I_N(f_{2^k M} g_M)-(I_Nf_{2^k M})g_M\|_r
&\lesssim 2^{-k} \bigl(\frac{2^kM}{N}\bigr)^{s-1}\|f_{2^k M}\|_{r_1}\|g_M\|_{r_2}\\
&\lesssim 2^{-k}M^{-(1-s+\nu)}\|I_Nf\|_{r_1}\||\nabla|^{1-s+\nu} g\|_{r_2}.
\end{align*}
Summing over $M$ and $k$ such that $N\lesssim 2^kM$, and recalling that $0<\nu<s$, we get
\begin{align}\label{d2}
III\lesssim N^{-(1-s+\nu)}\|I_Nf\|_{r_1}\||\nabla|^{1-s+\nu}g\|_{r_2}.
\end{align}

To estimate $I$, we apply the same argument as for $III$.  We get
\begin{align}\label{d3}
I=\|I_N(f_{\gtrsim N}g_{\le 1})-(I_Nf_{\gtrsim N})g_{\le 1} \|_r
&\lesssim \sum_{k\in \mathbb N,2^k \gtrsim N}\|I_N(f_{2^k}g_{\le 1})-(I_Nf_{2^k})g_{\le 1}\|_r \nonumber\\
&\lesssim \sum_{k\in \mathbb N,2^k\gtrsim N} 2^{-k}\|I_Nf\|_{r_1}\|g\|_{r_2}\nonumber\\
&\lesssim N^{-1}\|I_Nf\|_{r_1}\|g\|_{r_2}.
\end{align}
Putting \eqref{difference} through \eqref{d3} together, we derive \eqref{bilinear eq}.
\end{proof}

As an application of Lemma \ref{bilinear} we have the following commutator estimate:

\begin{lemma}\label{nablaif}
Let $1<r,r_1,r_2<\infty$ be such that $\frac 1r=\frac 1{r_1}+\frac 1{r_2}$.  Then, for any $0<\nu<s$ we have
\begin{align}
\|\nabla I_NF(u)-(I_N\nabla u)F'(u)\|_r\lesssim N^{-1+s-\nu}\|\nabla I_N u\|_{r_1}\|\langle\nabla\rangle^{1-s+\nu}F'(u)\|_{r_2}\label{na1}\\
\|\nabla I_N F(u)\|_r\lesssim \|\nabla I_N u\|_{r_1}\|F'(u)\|_{r_2}+N^{-1+s-\nu} \|\nabla I_Nu\|_{r_1}\|\langle\nabla\rangle^{1-s+\nu}F'(u)\|_{r_2}. \label{na2}
\end{align}
\end{lemma}
\begin{proof}
As
$$
\nabla F(u)=\nabla u \cdot F'(u),
$$
the estimate \eqref{na1} follows immediately from Lemma \ref{bilinear} with $f:=\nabla u$ and $g:=F'(u)$. The estimate \eqref{na2} is a
consequence of \eqref{na1} and the triangle inequality.
\end{proof}

Since we work at regularity $0<s<1$, we will need the following fractional chain rules to estimate our nonlinearity in $H_x^s$.
\begin{lemma}[Fractional chain rule for a $C^1$ function]\label{F C1}
Suppose that $F\in C^1(\mathbb C)$, $\sigma \in (0,1)$, and $1<r,r_1,r_2<\infty$ such that $\frac 1r=\frac 1{r_1}+\frac 1{r_2}$.   Then,
$$
\||\nabla|^{\sigma}F(u)\|_r\lesssim \|F'(u)\|_{r_1}\||\nabla|^{\sigma}u\|_{r_2}.
$$
\end{lemma}

\begin{lemma}[Fractional chain rule for a Lipschitz function]\label{F lip}
Let $F$ be a Lipschitz function, $\sigma \in (0,1)$, and $1<r<\infty$.  Then,
$$
\||\nabla|^{\sigma}F(u)\|_r\lesssim \|F'\|_{\infty}\||\nabla|^{\sigma}u\|_{r}.
$$
\end{lemma}

\begin{lemma}[Fractional derivatives for fractional powers]\label{fdfp}
Let $F$ be a H\"older continuous function of order $0<\alpha<1$.  Then, for every $0<\sigma<\alpha$, $1<r<\infty$,
and $\tfrac{\sigma}{\alpha}<\delta<1$ we have
\begin{align}\label{fdfp2}
\bigl\| |\nabla|^\sigma F(u)\bigr\|_r
\lesssim \bigl\||u|^{\alpha-\frac{\sigma}{\delta}}\bigr\|_{r_1} \bigl\||\nabla|^{\delta} u\bigr\|^{\frac{\sigma}{\delta}}_{\frac{\sigma}{\delta}r_2},
\end{align}
provided $\tfrac{1}{r}=\tfrac{1}{r_1} +\tfrac{1}{r_2}$ and $(1-\frac\sigma{\alpha {\delta}})r_1>1$.
\end{lemma}

The first two results originate in \cite{chris:weinstein} and \cite{staf}; for a textbook treatment see \cite{taylor}.  The third result
can be found in Appendix~A of \cite{MV}.
Using the chain rule estimates in these lemmas, we can upgrade the pointwise in time commutator estimate
in Lemma~\ref{nablaif} to a spacetime estimate:
\begin{lemma}\label{gif}
Let $I$ be a compact time interval and let $\frac 1{1+\min\{1,\frac 4d\}}<s<1$.  Then,
\begin{align}
\|(I_N\nabla u)F'(u)-\nabla I_N F(u)\|_{L_t^2L_x^{\frac{2d}{d+2}}(I\times \R^d)}
\lesssim N^{-\min\{1,\frac 4d\}s+}\|\lrn I_N u\|_{S^0(I)}^{1+\frac 4d}\label{gif2}\\
\|\lrn I_N F(u)\|_{N^0(I)}
\lesssim |I|^{\frac {2s}d}\|\lrn I_Nu\|_{S^0(I)}^{1+\frac 4d}+N^{-\min\{1,\frac 4d\}s+}\|\lrn I_N u\|_{S^0(I)}^{1+\frac 4d}. \label{gif1}
\end{align}
\end{lemma}

\begin{proof}
Throughout the proof, all spacetime norms will be taken on the slab $\ir$.

As by hypothesis $s>\frac 1{1+\min\{1,\frac 4d\}}$, there exists $\eps_0>0$ such that for any $0<\eps<\eps_0$ we have
$s>\frac {1+\eps}{1+\min\{1,\frac 4d\}}$.  Let $\nu:=\min\{1,\frac 4d\}s-(1-s)-\eps$ with $0<\eps<\eps_0$.
It is easy to check that $0<\nu<s$.  Applying Lemma~\ref{nablaif} with this value of $\nu$, we get
\begin{align*}
\|(\nabla I_N u)F'(u)-&\nabla I_N F(u)\|_{2,\frac {2d}{d+2}}\\
&\lesssim N^{-\min\{1,\frac 4d\}s+\eps}\|\nabla I_N u\|_{2,\frac {2d}{d-2}}\|\lrn^{\min\{1,\frac 4d\}s-\eps}F'(u)\|_{\infty,\frac d2}.
\end{align*}
The claim \eqref{gif2} will follow immediately from the estimate above, provided we show
\begin{align}
\|\langle \nabla \rangle^{\min\{1,\frac 4d\}s-\eps}F'(u)\|_{\infty,\frac d2}\lesssim \|\lrn I_N u\|_{\infty,2}^{\frac 4d}. \label{a4}
\end{align}

We start by observing that for any $(q,r)$ admissible pair,
\begin{align}
\|\langle \nabla\rangle^s u\|_{q,r} \lesssim \|\lrn I_N u\|_{q,r}.\label{a1}
\end{align}
Indeed, decomposing $u:=u_{\le N}+u_{>N}$ and using Lemma \ref{basic property} and the fact that $I_N$ is the identity operator on frequencies
$|\xi|\leq N$, we get
\begin{align*}
\|u\|_{q,r}
&\le \|u_{\le N}\|_{q,r}+\|u_{> N}\|_{q,r}\\
& \lesssim \|I_N u_{\le N}\|_{q,r}+N^{-1}\|\nabla I_N u_{>N}\|_{q,r}\\
&\lesssim \|\lrn I_N u\|_{q,r}.
\end{align*}
Similarly, we estimate
\begin{align*}
\||\nabla|^s u\|_{q,r}
&\le \||\nabla|^s u_{\le N}\|_{q,r} + \||\nabla|^s u_{>N}\|_{q,r}\\
&\lesssim \||\nabla|^s I_N u_{\le N}\|_{q,r} + N^{s-1}\|\nabla I_N u\|_{q,r}\\
&\lesssim \|\lrn I_Nu\|_{q,r},
\end{align*}
and the estimate \eqref{a1} follows.

As $F'(u)=O(|u|^{\frac 4d})$, by \eqref{a1} we get
\begin{align}
\|F'(u)\|_{\infty,\frac d2}\lesssim \|\lrn I_N u\|_{\infty,2}^{\frac 4d}\label{a2}.
\end{align}

We first prove \eqref{a4} for $d\leq 4$.  By \eqref{a2}, we estimate
\begin{align*}
\|\langle \nabla \rangle^{\min\{1,\frac 4d\}s-\eps}F'(u)\|_{\infty,\frac d2}
&=\|\langle \nabla \rangle^{s-\eps}F'(u)\|_{\infty,\frac d2}\\
&\lesssim \|\langle \nabla \rangle^{s}F'(u)\|_{\infty,\frac d2}\\
&\lesssim \|F'(u)\|_{\infty,\frac d2} + \||\nabla|^{s}F'(u)\|_{\infty,\frac d2}\\
&\lesssim \|\lrn I_N u\|_{\infty,2}^{\frac 4d}+ \||\nabla|^{s}F'(u)\|_{\infty,\frac d2}.
\end{align*}
Using Lemmas \ref{F C1} (for $d=3$) and \ref{F lip} (for $d=4$) together with \eqref{a1}, we estimate
\begin{align*}
\||\nabla|^sF'(u)\|_{\infty,\frac d2}
\lesssim \||\nabla|^s u\|_{\infty,2}\|u\|_{\infty,2}^{\frac{4}{d} -1}
\lesssim \|\lrn I_N u\|_{\infty,2}^{\frac 4d}.
\end{align*}
Thus, \eqref{a4} holds for $d\leq 4$.

If $d>4$, using Lemma~\ref{fdfp} (with $\alpha:=\frac 4d$, $\sigma:=\frac{4s}{d}-\eps$, and $\delta:=s$) and \eqref{a1}, we get
\begin{align*}
\||\nabla|^{\min\{1,\frac 4d\}s-\eps}F'(u)\|_{\infty,\frac d2}
=\||\nabla|^{\frac{4s}{d}-\eps}F'(u)\|_{\infty,\frac d2}
&\lesssim \|u\|_{\infty,2}^{\frac{\eps}{s}}\||\nabla|^s u\|_{\infty,2}^{\frac{4}{d}-\frac{\eps}{s}}\\
&\lesssim \|\lrn I_N u\|_{\infty,2}^{\frac{4}{d}}.
\end{align*}
From this and \eqref{a2} we derive \eqref{a4} in the case $d>4$.  The claim \eqref{gif2} follows.

We now consider \eqref{gif1}.  Using \eqref{estimate gamma-rho}, \eqref{i1}, and \eqref{a1}, we obtain
\begin{align*}
\|I_NF(u)\|_{N^0(I)}
\lesssim \|F(u)\|_{\gamma',\rho'}
\lesssim |I|^{\frac {2s}d}\|\lrs u\|_{\gamma,\rho}^{1+\frac 4d}
\lesssim |I|^{\frac {2s}d} \|\lrn I_N u\|_{\gamma,\rho}^{1+\frac 4d}.
\end{align*}
Similarly, by H\"older and \eqref{a1},
\begin{align*}
\|(I_N \nabla u)F'(u)\|_{N^0(I)}
&\lesssim \|(I_N \nabla u)F'(u)\|_{\gamma',\rho'}\\
&\lesssim |I|^{\frac{2s}{d}}\| I_N \nabla u\|_{\gamma, \rho}\|u\|_{\gamma, \rho^*}^{\frac 4d}\\
&\lesssim |I|^{\frac{2s}{d}}\|I_N \nabla u\|_{\gamma, \rho}\|\lrs u\|_{\gamma,\rho}^{\frac 4d}\\
&\lesssim |I|^{\frac{2s}{d}}\|\langle \nabla\rangle I_N u\|_{\gamma,\rho}^{1+\frac 4d}.
\end{align*}
The estimate \eqref{gif1} follows from the estimates above, \eqref{gif2} and the triangle inequality.
\end{proof}

We end this section with the following concentration-compactness lemma:
\begin{lemma}[Concentration-compactness, \cite{hk:revisited}]\label{compact}
Let $\{v_n\}_{n\geq 1}$ be a bounded sequence in $H^1(\R^d)$ such that
$$
\limsup_{n\to \infty}\|\nabla v_n\|_2\le M<\infty
$$
and
$$
\limsup_{n\to\infty}\|v_n\|_{2+\frac 4d}\ge m>0.
$$
Then, there exists $\{x_n\}_{n\geq 1}\subset \R^d$ such that, up to a subsequence,
$$
v_n(\cdot+x_n)\rightharpoonup  V \quad \text{weakly in } H^1_x \text{ as } n\to \infty
$$
with $\|V\|_2\ge \sqrt{\frac d{d+2}}\frac {m^{\frac 2d+1}}M\|Q\|_2$.
\end{lemma}

%
%
%
%

\section{local $H_x^s$ theory}\label{local}
In this section, we review the local $H_x^s$ theory for the equation \eqref{equation} as described in \cite{cazenave:book}.

\begin{proposition}[Local well-posedness in $H^s_x$, \cite{cazenave:book}]\label{Hs theory}
Let $ 0<s< 1$ and $u_0\in H^s(\R^d)$.  Then, the equation \eqref{equation} is wellposed on $[0,T_{lwp}]$ with
$$
T_{lwp}=C_0 \|\langle \nabla \rangle^su_0\|_2^{-\frac 2s}.
$$
Moreover, the unique solution $u$ enjoys the following estimate:
$$
\|\langle \nabla \rangle^s u\|_{S^0([0,T_{lwp}])}\lesssim \|\langle \nabla \rangle^s u_0\|_2.
$$
Here, $C_0$ and the implicit constant depend only on the dimension $d$ and the regularity $s$.
\end{proposition}

As a direct consequence of the $H_x^s$ local well-posedness result, we have the following lower bound on the blowup rate of the $H^s_x$-norm:
\begin{corollary}[Blowup criterion, \cite{cazenave:book}]\label{criterion}
Let $0<s< 1$ and $u_0\in H_x^s$.  Assume that the unique solution $u$ to \eqref{equation} blows up at time $0<T^*<\infty$.  Then, there exists a
constant $C$ depending only on $d$ and $s$ such that
$$
\|u(t)\|_{H_x^s}\ge C(T^*-t)^{-\frac s2}.
$$
\end{corollary}

As a variation of Proposition~\ref{Hs theory}, we have
\begin{proposition}[`Modified' local well-posedness]\label{mhs}
Let $\frac {\frac 4d+1}{1+\min\{1,\frac 4d\}+\frac 4d}<s<1$ and $u_0\in H_x^s$.  Let also
\begin{align}
N&\gg \|u_0\|_{H_x^s}^{\frac 4{\min\{4,d\}s-(4+d)(1-s)-d\eps}} \label{choose n} \quad \text{for any } \eps>0 \ \text{sufficiently small}\\
\tilde T_{lwp}&:= c_0\|\lrn I_N u_0\|_2^{-\frac 2s}\quad \mbox{for a small constant }  c_0=c_0(d,s). \label{choose t}
\end{align}
Then, \eqref{equation} is wellposed on $[0,\tilde T_{lwp}]$ and moreover
\begin{align}\label{s2}
\|\lrn I_N u\|_{S^0([0,\tilde T_{lwp}])}\lesssim \|\lrn I_N u_0\|_2.
\end{align}
\end{proposition}

\begin{proof}
As by Lemma \ref{basic property},
$$
\|u_0\|_{H^s_x}\lesssim \|\lrn I_N u_0\|_2,
$$
choosing $c_0$ sufficiently small, we get
$$
\tilde T_{lwp}\le T_{lwp}.
$$
Thus, \eqref{equation} is wellposed in $H_x^s$ on $[0,\tilde T_{lwp}]$. Let $u$ be the unique solution; by Proposition~\ref{Hs theory}
we have
\begin{equation}
\|\lrs u\|_{S^0([0,\tilde T_{lwp}])} \lesssim \|u_0\|_{H_x^s}. \label{u estimate}
\end{equation}

On the other hand, by Strichartz, for any $t\le \tilde T_{lwp}$ we estimate
$$
\|\lrn I_N u\|_{S^0([0,t])}\lesssim \|\lrn I_Nu_0\|_2+\|\lrn I_N F(u)\|_{N^0([0,t])}.
$$
As by hypothesis, $s>\frac {\frac 4d+1}{1+\min\{1,\frac 4d\}+\frac 4d}>\frac 1{1+\min\{1,\frac 4d\}}$, by Lemma~\ref{gif} we get
\begin{align}
\|\lrn I_Nu\|_{S^0([0,t])}
&\lesssim \|\lrn I_N u_0\|_2+t^{\frac {2s}d}\|\lrn I_N u\|_{S^0([0,t])}^{1+\frac 4d} \nonumber\\
& \quad +N^{-\min\{1,\frac 4d\}s+\eps} \|\lrn I_N u\|_{S^0([0,t])}^{1+\frac 4d}\label{s1}
\end{align}
for any $\eps>0$ sufficiently small.
For $N$ satisfying \eqref{choose n}, we use Lemma \ref{basic property} and \eqref{u estimate} to estimate
\begin{align*}
N^{-\min\{1,\frac 4d\}s+\eps}\|& \lrn I_N u\|_{S^0([0,t])}^{1+\frac 4d}\\
&\lesssim N^{-\min\{1,\frac 4d\}s+\eps}N^{(1+\frac 4d)(1-s)}\|\lrs u\|_{S^0([0,t])}^{1+\frac 4d}\\
&\ll \|u_0\|_{H_x^s}^{\frac 4{\min\{4,d\}s-(4+d)(1-s)-d\eps}\times [-\min\{1,\frac 4d\}s+\eps+(1+\frac 4d)(1-s)]}\|u_0\|_{H_x^s}^{1+\frac 4d} \\
&\lesssim \|u_0\|_{H_x^s}^{-\frac 4d}\|u_0\|_{H_x^s}^{1+\frac 4d}
\lesssim \|\lrn I_N u_0\|_2.
\end{align*}
Here, we have used the fact that $s>\frac {\frac 4d+1}{1+\min\{1,\frac 4d\}+\frac 4d}$ implies that there exists $\eps_1>0$ sufficiently small such that
for any $0<\eps<\eps_1$ we have $s>\frac {\frac 4d+1+\eps}{1+\min\{1,\frac 4d\}+\frac 4d}$; thus the power of $N$ in the estimates above
is negative for $0<\eps<\eps_1$.

Returning to \eqref{s1}, we conclude that
$$
\|\lrn I_N u\|_{S^0([0,t])}
\lesssim \|\lrn I_N u_0\|_2+t^{\frac {2s}d}\|\lrn I_N u\|_{S^0([0,t])}^{1+\frac 4d}.
$$
Standard arguments yield \eqref{s2}, provided $$t\leq \tilde T_{lwp}=c_0(d,s)\|\lrn I_N u_0\|_2^{-\frac 2s}.$$
\end{proof}

%
%
%
%

\section{Modified energy increment}
The main purpose of this section is to prove that the modified energy of $u$, $E(I_Nu)$, grows much slower than the modified kinetic
energy of $u$, $\|\nabla I_N u\|_2^2$. As will be shown later, this result is crucial in establishing the main theorems.

Before stating the result, we need to introduce more notation.  We define
$$
 \Lambda(t):=\sup_{0\le \tau\le t} \|u(\tau)\|_{H_x^s} \quad \text{and} \quad \Sigma(t):=\sup_{0\le \tau\le t} \|I_N u(\tau)\|_{H^1_x}.
$$

With this notation we have the following
\begin{proposition}[Increment of the modified energy]\label{slower}
Let $s_0(d)<s< 1$ and let $u_0\in H_x^s$ such that the corresponding solution $u$ to \eqref{equation} blows up at time $0<T^*<\infty$.
Let $0<T<T^*$.  Then, for
\begin{align}\label{N(T)}
N(T):=C \Lambda(T)^{\frac {p(s)}{2(1-s)}},
\end{align}
we have
$$
|E(I_{N(T)}u(T))|\lesssim \Lambda(T)^{p(s)}.
$$
Here, $C$ and the implicit constant depend only on $s$, $T^*$, and $\|u_0\|_{H_x^s}$, and $p(s)$ is given by
$$
p(s):=\frac {2(2+\frac 2s+\frac 8d)(1-s)}{\min\{1,\frac 4d\}s-(\frac 2s+\frac 8d)(1-s)-}.
$$
\end{proposition}

Note that by Lemma~\ref{basic property}, $\Lambda(T)\lesssim \Sigma(T)$.  Thus, if the solution blows up at time $T^*$, the modified energy
$E(I_{N(T)}u(T))$ is at most $O(\Sigma(T)^{p(s)})$, which is much smaller than the modified kinetic energy,
$\|\nabla I_{N(T)} u(T)\|_{2}^2=O(\Sigma(T)^2)$ for $s>s_0(d)$.

We prove Proposition \ref{slower} in two steps.  The first step is to control the increment of the modified energy of $u$ on intervals of local
well-posedness $[0,\tilde T_{lwp}]$.  The second step is to divide the interval $[0,T]$ into finitely many subintervals of local well-posedness,
control the increment of the modified energy of $u$ on each of these subintervals, and sum these bounds.

We start with the following
\begin{lemma}[Local increment of the modified energy]\label{increment}
Let $\frac {1+\frac 4d}{1+\frac 4d+\min\{1,\frac 4d\}}<s< 1$ and let $u_0\in H_x^s$. Assume that $N$ and $\tilde T_{lwp}$ satisfy \eqref{choose n}
and \eqref{choose t} respectively, that is,
\begin{align*}
N&\gg \|u_0\|_{H_x^s}^{\frac 4{\min\{4,d\}s-(4+d)(1-s)-}}\\
\tilde T_{lwp}&:= c_0\|\lrn I_N u_0\|_2^{-\frac 2s}\quad \mbox{for a small constant }  c_0=c_0(d,s).
\end{align*}
Then,
$$
\sup_{t\in [0,\tilde T_{lwp}]}|E(I_Nu(t))|\leq |E(I_N u_0)|+ CN^{-\min\{1,\frac 4d\}s+}(\|I_N\lrn u_0\|_2^{2+\frac4d}+\|I_N\lrn u_0\|_2^{2+\frac 8d}).
$$
Here, the constant $C$ depends on $s$, $T^*$, and $\|u_0\|_{H_x^s}$.
\end{lemma}

\begin{proof}
Note that by Proposition~\ref{mhs}, \eqref{equation} is wellposed on $[0,\tilde T_{lwp}]$.  Furthermore, the unique solution $u$ to \eqref{equation}
on $[0,\tilde T_{lwp}]$ satisfies
\begin{equation}
\|\lrn I_N u\|_{S^0([0,\tilde T_{lwp}])}\lesssim \|\lrn I_N u_0\|_2. \label{e1}
\end{equation}

Let $0<t\le \tilde T_{lwp}$; throughout the rest of the proof, all spacetime norms will be taken on $[0,t]\times\R^d$.
By the Fundamental Theorem of Calculus, we can write the modified energy increment as
\begin{align*}
E(I_Nu(t))-E(I_Nu_0)&=\int_0^t\frac {\partial}{\partial s}E(I_Nu(s))\,ds\\
&=Re \int_0^t\int_{\R^d} \overline{I_N u_t}(-\Delta I_N u+F(I_N u))\,dx\,ds.
\end{align*}
As $I_N u_t=i\Delta I_N u-iI_N F(u)$, we have
$$
Re\int_0^t\int_{\R^d}\overline {I_Nu_t}(-\Delta I_N u+I_NF(u))\,dx\,ds=0.
$$
Thus, after an integration by parts,
\begin{align}
E(I_N u(t))-E(I_Nu_0)
&=Re \int_0^t\int_{\R^d}\overline {I_Nu_t}[F(I_N u)-I_NF(u)]\,dx\,ds\notag\\
&=-Im \int_0^t\int_{\R^d} \overline{\nabla I_Nu}\cdot\nabla[F(I_Nu)-I_NF(u)]\,dx\,ds\label{from lap}\\
&\quad -Im \int_0^t\int_{\R^d} \overline{I_NF(u)}\cdot [F(I_Nu)-I_NF(u)]\,dx\,ds \label{from nol}.
\end{align}

Consider the contribution from \eqref{from lap}.  By the triangle inequality,
\begin{align*}
\|\nabla[F(I_Nu)-I_NF(u)]\|_{2,\frac{2d}{d+2}}
&\lesssim \|(\nabla I_N u)[F'(I_N u)-F'(u)]\|_{2,\frac{2d}{d+2}}\\
&\quad +\|(\nabla I_N u) F'(u)-\nabla I_N F(u)\|_{2,\frac{2d}{d+2}}.
\end{align*}
By H\"older, \eqref{holder continuity}, \eqref{i2}, and \eqref{a1}, we estimate
\begin{align*}
\|(\nabla I_N u)&[F'(I_N u)-F'(u)]\|_{2,\frac{2d}{d+2}}\\
&\lesssim \|\nabla I_N u\|_{2,\frac{2d}{d-2}}\|F'(I_N u)-F'(u)\|_{\infty,\frac d2}\\
&\lesssim \|\lrn I_N u\|_{S^0([0,t])} \bigl\||I_Nu-u|^{\min\{1,\frac 4d\}}(|I_Nu|+|u|)^{\frac 4d-\min\{1,\frac 4d\}}\bigr\|_{\infty,\frac d2} \\
&\lesssim \|\lrn I_N u\|_{S^0([0,t])}\|P_{>N}u\|_{\infty,2}^{\min\{1,\frac 4d\}}\|u\|_{\infty,2}^{\frac 4d-\min\{1,\frac 4d\}}\\
&\lesssim N^{-\min\{1,\frac 4d\}}\|\lrn I_N u\|_{S^0([0,t])}^{1+\frac 4d}.
\end{align*}
Combining this with \eqref{gif2}, we get
\begin{align}\label{dif energy incr}
\|\nabla[F(I_Nu)-I_NF(u)]\|_{2,\frac{2d}{d+2}}&\lesssim N^{-\min\{1,\frac 4d\}s+}\|\lrn I_N u\|_{S^0([0,t])}^{1+\frac 4d}.
\end{align}

Therefore
\begin{align}
|\ref{from lap}|
&\lesssim \|\nabla I_N u\|_{2,\frac {2d}{d-2}}\|\nabla[F(I_Nu)-I_NF(u)]\|_{2,\frac{2d}{d+2}}\notag\\
&\lesssim N^{-\min\{1,\frac 4d\}s+}\|\lrn I_N u\|_{S^0([0,t])}^{2+\frac 4d}.\label{fl1}
\end{align}

We turn now toward \eqref{from nol}.  By \eqref{dif energy incr} and Sobolev embedding, we estimate
\begin{align}
|\eqref{from nol}|
&\lesssim \|\nabla[F(I_Nu)-I_NF(u)]\|_{2,\frac{2d}{d+2}}\||\nabla|^{-1}I_NF(u)\|_{2,\frac {2d}{d-2}}\nonumber \\
&\lesssim  N^{-\min\{1,\frac 4d\}s+}\|\lrn I_N u\|_{S^0([0,t])}^{1+\frac 4d}\|I_NF(u)\|_{2,2}.\label{fn1}
\end{align}
To estimate the last factor in \eqref{fn1}, we drop the operator $I_N$ and use Sobolev embedding to obtain
$$
\|I_NF(u)\|_{2,2}
\lesssim \|u\|_{\frac {2(d+4)}d, \frac{2(d+4)}d}^{1+\frac 4d}
\lesssim \||\nabla|^{\frac d{d+4}}u\|_{\frac {2(d+4)}d, \frac{2(d+4)}{d+2}}^{1+\frac 4d}.
$$
Note that $(\tfrac{2(d+4)}d,\tfrac{2(d+4)}{d+2})$ is a Schr\"odinger
admissible pair.  Decompose $u:=u_{\leq N}+u_{>N}$. To estimate the
low frequencies, we use the fact that $I_N$ is the identity on
frequencies $|\xi|\leq N$:
$$
\||\nabla|^{\frac d{d+4}}u_{\le N}\|_{\frac {2(d+4)}d,\frac{2(d+4)}{d+2}}\lesssim \|\lnr I_Nu\|_{S^0([0,t])}.
$$
For the high frequencies, we use Lemma~\ref{basic property} to get
$$
\||\nabla| ^{\frac d{d+4}}u_{> N}\|_{\frac{2(d+4)}d,\frac{2(d+4)}{d+2}}
\lesssim N^{-\frac 4{d+4}} \|\nabla I_Nu_{>N}\|_{S^0([0,t])},
$$
provided $s>\frac d{d+4}$; this condition is satisfied since by assumption,
$$
s>\tfrac {1+\frac 4d}{1+\frac 4d+\min\{1,\frac 4d\}}>\tfrac d{d+4}.
$$
Therefore,
\begin{equation}
\|I_NF(u)\|_{2,2}\lesssim \|\lrn I_Nu\|_{S^0([0,t])}^{1+\frac 4d}.\label{fn3}
\end{equation}
By \eqref{fn1} and \eqref{fn3}, we obtain
\begin{align}
|\eqref{from nol}|&\lesssim N^{-\min\{1,\frac 4d\}s+}\|\lrn I_N u\|_{S^0([0,t])}^{2+\frac 8d}. \label{fn4}
\end{align}

Collecting \eqref{e1}, \eqref{fl1}, and \eqref{fn4}, we get
$$
|E(I_Nu(t))-E(I_Nu_0)|
\lesssim N^{-\min\{1,\frac 4d\}s+}(\|\lrn I_N u_0\|_2^{2+\frac 4d}+\|\lrn I_N u_0\|_2^{ 2+\frac 8d}).
$$
This proves Lemma~\ref{increment}.
\end{proof}

Next, we use Lemma~\ref{increment} to prove Proposition \ref{slower}.

Let $T<T^*$ and $\Lambda(T)$ and $\Sigma(T)$ defined as in the beginning of this section. By Proposition~\ref{mhs}, if we take
\begin{equation}\label{sl1}
\begin{cases}
N(T)\gg \Lambda(T)^{\frac 4{\min\{4,d\}s-(4+d)(1-s)-}} \\
\delta:=c_0\Sigma(T)^{-\frac 2s},
\end{cases}
\end{equation}
then the solution $u$ satisfies the estimate
$$
\|\lnr I_{N(T)}u\|_{S^0([t, t+\delta])}\lesssim \|\lrn I_{N(T)}u(t)\|_2\lesssim \Sigma(T),
$$
uniformly in $t$, provided $[t,t+\delta]\subset [0,T]$. Thus, splitting $[0,T]$ into $O(\frac T{\delta})$ subintervals and applying
Lemma~\ref{increment} on each of these subintervals, we get
\begin{align}
\sup_{t\in[0,T]}|E(I_{N(T)}u(t))|
&\lesssim |E(I_{N(T)}u_0)|+\frac T{\delta} N(T)^{-\min\{1,\frac 4d\}s+} \Sigma(T)^{2+\frac 4d}\notag\\
&\quad +\frac T{\delta}N(T)^{-\min\{1,\frac 4d\}s+}\Sigma(T)^{2+\frac 8d}\notag\\
&\lesssim |E(I_{N(T)}u_0)|+ N(T)^{-\min\{1,\frac 4d\}s+} \Sigma(T)^{2+\frac 4d+\frac 2s}\notag\\
&\quad + N(T)^{-\min\{1,\frac 4d\}s+}\Sigma(T)^{2+\frac 8d+ \frac 2s}.\label{increase}
\end{align}
Using interpolation, Sobolev embedding, and Lemma~\ref{basic property}, we estimate
\begin{align}
|E(I_{N(T)}u_0)|
&\lesssim \|\nabla I_N u_0\|_2^2 +\|I_N u_0\|_{2+\frac 4d}^{2+\frac 4d}\notag\\
&\lesssim N^{2(1-s)}\|u_0\|_{H^s_x}^2 + \|I_N u_0\|_2^{\frac 4d} \|\nabla I_N u_0\|_2^2\notag\\
&\lesssim N^{2(1-s)}\bigl(\|u_0\|_{H^s_x}^2 +\|u_0\|_{H^s_x}^{2+\frac 4d}\bigr)\notag\\
&\lesssim N^{2(1-s)}.\label{initial data}
\end{align}
Moreover, by Lemma~\ref{basic property}, we also have
\begin{align}\label{sigma lambda}
\Sigma(T)\lesssim N(T)^{1-s}\Lambda(T).
\end{align}
Substituting \eqref{initial data} and \eqref{sigma lambda} into \eqref{increase}, we obtain
\begin{align}
\sup_{t\in[0,T]}|E(I_{N(T)}u(t))|
&\lesssim N(T)^{2(1-s)}+N(T)^{-\min\{1,\frac 4d\}s+(2+\frac 4d+\frac 2s)(1-s)+}\Lambda(T)^{2+\frac 4d+\frac 2s}\nonumber\\
&\quad +N(T)^{-\min\{1,\frac 4d\}s+(2+\frac 8d+\frac 2s)(1-s)+}\Lambda(T)^{2+\frac 8d+\frac 2s}. \label{sl2}
\end{align}
Optimizing \eqref{sl2}, we observe that if
\begin{equation}
N(T)\sim\Lambda(T)^{\frac {2+\frac 8d+\frac 2s}{\min\{1,\frac 4d\}s-(1-s)(\frac 8d+\frac 2s)-}}, \label{sl3}
\end{equation}
then $N(T)$ satisfies the assumption \eqref{sl1} and moreover,
\begin{align*}
\sup_{t\in[0,T]}|E(I_{N(T)}u(t))|&\lesssim N(T)^{2(1-s)}\\
&\lesssim \Lambda(T)^{\frac {2(2+\frac 8d+\frac 2s)(1-s)}{\min\{1,\frac 4d\}s-(1-s)(\frac 8d+\frac 2s)-}}.
\end{align*}
Let
$$
p(s):=\frac {2(2+\frac 8d+\frac 2s)(1-s)}{\min\{1,\frac 4d\}s-(1-s)(\frac 8d+\frac 2s)-}.
$$
Then, a little work shows that the condition $0<p(s)<2$ leads to the restriction
$$
s>s_0(d).
$$
Thus, for $N(T)$ defined in \eqref{sl3} and $s>s_0(d)$, we have $0<p(s)<2$ and
$$
\sup_{t\in [0,T]}|E(I_{N(T)}u)(t)|\lesssim \Lambda(T)^{p(s)}.
$$
This proves Proposition~\ref{slower}.

%
%
%
%

\section{Proof of Theorem~\ref{blowup}}
In this section, we use Proposition~\ref{slower} together with Lemma~\ref{compact} to prove Theorem~\ref{blowup}.

We choose a sequence of times $\{t_n\}_{n\geq 1}$, such that $t_n\to T^*$ as $n\to\infty$ and
$$
\|u(t_n)\|_{H_x^s}=\Lambda(t_n).
$$
As the solution $u$ blows up at time $T^*$, we must have $\Lambda(t_n)\to \infty$ as $n\to \infty$.

Set
$$
\psi_n(x):=\rho_n^{\frac d2}(I_{N(t_n)}u)(t_n,\rho_nx),
$$
where $N(t_n)$ is given by \eqref{N(T)} with $T:=t_n$ and the parameter $\rho_n$ is given by
$$
\rho_n:=\frac {\|\nabla Q\|_2}{\|\nabla I_{N(t_n)}u(t_n)\|_2}.
$$
By Lemma~\ref{basic property} and Corollary~\ref{criterion}, we get
$$
\rho_n\lesssim \frac {1}{\|u(t_n)\|_{H_x^s}}\lesssim (T^* -t_n)^{\frac s2}.
$$
Basic calculations show that $\{\psi_n\}_{n\geq 1}$ is a bounded sequence in $H^1_x$.  Indeed,
\begin{align}
\|\psi_n\|_2&= \|I_{N(t_n)}u(t_n)\|_2\le \|u(t_n)\|_2=\|u_0\|_2\nonumber\\
\|\nabla \psi_n\|_2&=\rho_n\|I_{N(t_n)}\nabla u(t_n)\|_2=\|\nabla Q\|_2.\label{kinetic}
\end{align}

By Proposition~\ref{slower} (with $T=t_n$), we can estimate the energy of $\psi_n$ as follows:
\begin{align*}
E(\psi_n)=\rho_n^2E(I_{N(t_n)}u(t_n))
\lesssim\rho_n^2\Lambda(t_n)^{p(s)}\lesssim \|u(t_n)\|_{H_x^s}^{p(s)-2}.
\end{align*}
Thus, as $p(s)<2$ for $s>s_0(d)$,
$$
E(\psi_n)\to 0 \quad \text{as} \quad n\to \infty,
$$
which by the definition of the energy and \eqref{kinetic} implies
\begin{equation}
\|\psi_n\|_{2+\frac 4d}^{2+\frac 4d}\to \frac {d+2}d\|\nabla Q\|_2^2\ \ \mbox{as  } n\to \infty.\label{potential limit}
\end{equation}

Applying Lemma~\ref{compact} to the sequence $\{\psi_n\}_{n\geq 1}$ (with $M:=\|\nabla Q\|_2$
and $m:=(\frac {d+2}d\|\nabla Q\|_2^2)^{\frac d{2d+4}}$), we derive the existence of a sequence $\{x_n\}_{n\geq 1}\subset \R^d$ and of a function
$V\in H^1(\R^d)$ such that $\|V\|_2\ge \|Q\|_2$ and, up to a subsequence,
$$
\psi_n(\cdot+x_n)\rightharpoonup V \quad \text{weakly in } H^1_x \quad \text{as } n\to\infty,
$$
that is,
\begin{align}\label{limit1}
\rho_n^{\frac d2}(I_{N(t_n)}u)(t_n,\rho_n\cdot +x_n)\rightharpoonup V \quad \text{weakly in } H^1_x \quad \text{as } n\to\infty.
\end{align}

To prove Theorem~\ref{blowup}, we have to eliminate the smoothing operator $I_{N(t_n)}$ from \eqref{limit1}.
We do so at the expense of trading the weak convergence in $H^1_x$ for convergence in the sense of distributions.  Indeed,
for any $\sigma<s$ we have
\begin{align}
\|\rho_n^{\frac d2}\bigl(u(t_n)-I_{N(t_n)}u(t_n)\bigr)&(\rho_n\cdot +x_n)\|_{\dot H_x^{\sigma}}\notag\\
&=\rho_n^{\sigma}\|P_{\ge N(t_n)}u(t_n)\|_{\dot H_x^{\sigma}}\nonumber\\
&\lesssim \rho_n^{\sigma}N(t_n)^{\sigma-s}\|P_{\ge N(t_n)}u(t_n)\|_{\dot H_x^s}\nonumber\\
&\lesssim \Lambda(t_n)^{-\sigma}\Lambda(t_n)^{\frac {(\sigma-s)p(s)}{2(1-s)}}\|P_{\ge N(t_n)}u(t_n)\|_{ H_x^s}\nonumber\\
&\lesssim \Lambda(t_n)^{1-\sigma+\frac {(\sigma-s)p(s)}{2(1-s)}}.\label{exponent}
\end{align}
Plugging the explicit expression for $p(s)$ in the above computation, we find that for
$$
\sigma<\tilde s:=\frac {2d+8s+s^2d(2-\min\{1,\frac 4d\})}{4d+16s-s^2(d\min\{1,\frac 4d\}+8)},
$$
the exponent of $\Lambda(t_n)$ in \eqref{exponent} is negative.  Hence,
\begin{equation}
\|\rho_n^{\frac d2}\bigl(u(t_n)-I_{N(t_n)}u(t_n)\bigr)(\rho_n \cdot+x_n)\|_{H^{\tilde s-}_x} \to 0 \ \mbox{ as } n\to \infty.\label{limit2}
\end{equation}
Combining \eqref{limit1} and \eqref{limit2} finishes the proof of Theorem~\ref{blowup}.

%
%
%
%

\section{Proof of Theorem~\ref{concentration}}
By Theorem~\ref{blowup}, there exists a blowup profile $V\in H_x^1$, with
$\|V\|_2\ge\|Q\|_2$, and there exist sequences $\{t_n,\rho_n,x_n\}_{n\geq 1}\subset \R_+\times\R_+^*\times \R^d$ such that $t_n\to T^*$,
\begin{equation}
\frac {\rho_n}{(T^*-t_n)^{\frac s2}} \lesssim 1 \quad \text{for all } n\geq 1, \label{rhon}
\end{equation}
and
\begin{align}\label{limit4}
\rho_n^{\frac d2}u(t_n,\rho_n \cdot+x_n)\rightharpoonup V \quad \text{weakly as } n\to\infty.
\end{align}

From \eqref{limit4} it follows that for any $R>0$ we have
$$
\liminf_{n\to \infty} \rho_n^d\int_{|x|\le R}|u(t_n,\rho_n x+x_n)|^2\ge \int_{|x|\le R}|V|^2dx,
$$
which, by a change of variables, yields
$$
\liminf_{n\to \infty} \sup_{y\in \R^d}\int_{|x-y|\le R\rho_n}|u(t_n,x)|^2dx\ge \int_{|x|\le R}|V|^2dx.
$$
As by hypothesis $\frac {(T^*-t_n)^{\frac s2}}{\alpha(t_n)}\to 0$ as $n\to \infty$, \eqref{rhon} implies that $\frac {\rho_n}{\alpha(t_n)}\to 0$
as $n\to \infty$.  Therefore,
$$
\liminf_{n\to \infty} \sup_{y\in \R^d}\int_{|x-y|\le \alpha(t_n)}|u(t_n,x)|^2dx\ge \int_{|x|\le R}|V|^2dx.
$$
Letting $R\to \infty$, we obtain
$$
\liminf_{n\to \infty} \sup_{y\in \R^d}\int_{|x-y|\le \alpha(t_n)}|u(t_n,x)|^2dx\ge \|V\|_2^2.
$$
As $\|V\|_2\geq \|Q\|_2$, this implies
$$
\limsup_{t\to T^*} \sup_{y\in \R^d}\int_{|x-y|\le \alpha(t)}|u(t,x)|^2dx\ge \|Q\|_2^2.
$$
As for any fixed time $t$, the map $y\to \int_{|x-y|\le \alpha(t)}|u(t,x)|^2dx$ is continuous and goes to zero as $y\to \infty$,
there exists $x(t)\in \R^d$ such that
$$
\sup_{y\in\R^d}\int_{|x-y|\le \alpha(t)}|u(t,x)|^2dx=\int_{|x-x(t)|\le\alpha(t)} |u(t,x)|^2dx.
$$
This finally implies
$$
\limsup_{t\to T^*}\int_{|x-x(t)|\le \alpha(t)}|u(t,x)|^2dx\ge \|Q\|_2^2,
$$
which proves Theorem~\ref{concentration}.

%
%
%
%

\section{Proof of Theorem~\ref{gs blowup}}
In this section, we upgrade Theorem~\ref{blowup} to Theorem~\ref{gs blowup} under the additional assumption $\|u_0\|_2=\|Q\|_2$.

With the notation used in the proof of Theorem~\ref{blowup}, we have
$$
\|\psi_n\|_2\le \|u_0\|_2=\|Q\|_2\le \|V\|_2.
$$
On the other hand, using the semi-continuity of weak convergence,
$$
\|V\|_2\le \liminf_{n\to \infty}\|\psi_n\|_2\le\|Q\|_2.
$$
Therefore,
$$
\|V\|_2=\|Q\|_2=\lim_{n\to \infty} \|\psi_n\|_2.
$$
Thus, as $\psi_n(\cdot+x_n) \rightharpoonup V$ weakly in $L_x^2$ (up to a subsequence which we still denote by $\psi_n(\cdot+x_n)$), we conclude that
$$
\psi_n(\cdot+x_n)\to V \ \mbox{ strongly in } L_x^2.
$$

Moreover, by the Gagliardo-Nirenberg inequality and the boundedness of $\{\psi_n\}_{n\geq 1}$ in $H_x^1$, we have
$$
\psi_n(\cdot+x_n)\to V \ \mbox{ in } L_x^{2+\frac 4d}.
$$
Combining this with \eqref{potential limit} and the sharp Gagliardo-Nirenberg inequality, we obtain
$$
\|\nabla Q\|_2\le \|\nabla V\|_2.
$$
By the semi-continuity of weak convergence, we also have
$$
\|\nabla V\|_2\le \liminf_{n\to\infty}\|\nabla \psi_n\|_2=\|\nabla Q\|_2,
$$
and so
$$
\|\nabla V\|_2=\|\nabla Q\|_2= \lim_{n\to \infty}\|\nabla \psi_n\|_2.
$$
Thus, as $\psi_n(\cdot+x_n) \rightharpoonup V$ in $H^1_x$, we conclude that
$$
\psi_n(\cdot+x_n)\to V \ \mbox{ strongly in } H^1_x.
$$
In particular, this implies
$$
E(V)=0.
$$

Collecting the properties of $V$ we find
$$
V\in H_x^1,\quad \|V\|_2=\|Q\|_2,\quad \|\nabla V\|_2=\|\nabla Q\|_2, \quad \text{and}\quad E(V)=0.
$$
The variational characterization of the ground state, \cite{weinstein:ground state}, implies that
$$
V(x)=e^{i\theta} Q(x+x_0)
$$
for some $(e^{i\theta},x_0)\in (S^1\times \R^d)$.
Thus,
\begin{align}\label{limit3}
\rho_n^{\frac d2}(I_{N(t_n)}u)(t_n,\rho_nx+x_n)\to e^{i\theta}Q(x+x_0) \quad \text{strongly in }H^1_x \text{ as } n\to \infty.
\end{align}

Theorem~\ref{gs blowup} follows from \eqref{limit2} and \eqref{limit3}.


\begin{thebibliography}{10}

\bibitem{begout:vargas}
P.~Begout, A.~Vargas, \emph{Mass concentration phenomena for the $L^2$-critical nonlinear Schr\"odinger equation}, preprint.

\bibitem{bourgain:blowup}
J.~Bourgain, \emph{Refinements of Strichartz inequality and applications to 2D-NLS with critical nonlinearity}, I.M.R.N. \textbf{5} (1998), 253-283.

\bibitem{Blowscatter}
J. Bourgain, \emph{Scattering in the energy space and below for 3D NLS}, J. d'Analyse Math. \textbf{75} (1998), 267-297.

\bibitem{borg:book}
J. Bourgain, \emph{New global well-posedness results for non-linear Schr\"odinger equations}, AMS
Publications (1999).

\bibitem{cwI}
T. Cazenave, F.B. Weissler, \emph{Critical nonlinear Schr\"odinger Equation}, Non. Anal. TMA \textbf{14}
(1990), 807--836.

\bibitem{cazenave:book}
T. Cazenave, \textit{Semilinear Schr\"odinger equations,} Courant Lecture Notes in Mathematics, \textbf{10},
American Mathematical Society, 2003.

\bibitem{chris:weinstein} M. Christ, M. Weinstein, \emph{Dispersion of small amplitude solutions of the
generalized Korteweg-de Vries equation}, J. Funct. Anal. \textbf{100} (1991), 87-109.

\bibitem{coifmey:1}
R. R. Coifman, Y. Meyer, \emph{On commutators of singular integrals
and bilinear singular integrals}, Trans. AMS \textbf{212} (1975),
315--331.

\bibitem{coifmey:2}
R. R. Coifman, Y. Meyer, \emph{Ondelettes and operateurs III,
Operateurs multilineaires}, Actualites Mathematiques, Hermann, Paris
(1991).

\bibitem{ckstt:low1}
J. Colliander, M. Keel, G. Staffilani, H. Takaoka, T. Tao, \emph{Almost conservation laws and global
rough solutions to a nonlinear Schr\"odinger equation}, Math. Res. Lett. \textbf{9} (2002), 659--682.

\bibitem{ckstt:low4}
J. Colliander, M. Keel, G. Staffilani, H. Takaoka, T. Tao, \emph{Polynomial upper bounds for the orbital instability of the 1D cubic NLS below
the energy norm}, Discrete Contin. Dyn. Syst. \textbf{9} (2003), 31--54.

\bibitem{ckstt:low7}
J. Colliander, M. Keel, G. Staffilani, H. Takaoka, T. Tao, \emph{Polynomial upper bounds for the instability
of the nonlinear Schr\"odinger equation below the energy norm}, Commun. Pure Appl. Anal. \textbf{2} (2003), 33--50.

\bibitem{crsw}
J.~Colliander, S.~Raynor, C.~Sulem, J. D.~Wright, \emph{Ground state mass concentration in the $L\sp 2$-critical nonlinear Schr\"odinger
equation below $H\sp 1$}, Math. Res. Lett. \textbf{12} (2005), 357--375.

\bibitem{sijia}
D. Fang, S. Zhong, \emph{Cauchy problem for the $L^2$-critical nonlinear Schr\"odinger equation below $H\sp 1$},
Nonlinear Anal. \textbf{62} (2005), 117--130.

\bibitem{gerard}
P.~Gerard, \emph{Description du defaut de compacite de l'injection de Sobolev}, ESAIM Control Optim. Calc. Var. \textbf{3} (1998), 213--233.

\bibitem{gv:scatter}
J. Ginibre, G. Velo, \emph{Scattering theory in the energy space for a class of nonlinear Schr\"odinger
equations}, J. Math. Pure. Appl. \textbf{64} (1985), 363--401.

\bibitem{gv:Strichartz}
J. Ginibre, G. Velo, \emph{Smoothing properties and retarded estimates for some dispersive evolution equations},
Comm. Math. Phys. \textbf{144} (1992), 163--188.

\bibitem{glassey}
R. T. Glassey, \emph{On the blowing up of solution to the Cauchy problem for nonlinear Schr\"odinger
operators}, J. Math. Phys. \textbf{8} (1977), 1794--1797.

\bibitem{hayashi:tsutsumi}
N. Hayashi, Y. Tsutsumi, \emph{Remarks on the scattering problem for nonlinear Schr\"odinger equations},
Diff. Eq. Math. Phys. (1986), 162--168.

\bibitem{hk:revisited}
T.~Hmidi, S.~Keraani, \emph{Blowup theory for the critical nonlinear Schr\"odinger equations revisited},
Int. Math. Res. Not. \textbf{46} (2005), 2815--2828.

\bibitem{hk}
T.~Hmidi, S.~Keraani, \emph{Remarks on the blowup for the $L^2$-critical nonlinear Schr\"odinger equations}, preprint.

\bibitem{kato}
T. Kato, \emph{On nonlinear Schr\"odinger equations},  Ann. Inst. H. Poincare Phys. Theor.
\textbf{46}  (1987),  113--129.

\bibitem{katounique}
T. Kato, \emph{On nonlinear Schr\"odinger equations, II.  $H^s$-solutions and unconditional well-posedness}, J. d'Analyse Math. \textbf{67} (1995), 281--306.

\bibitem{tao:keel}
M. Keel, T. Tao, \emph{Endpoint Strichartz Estimates}, Amer. Math. J. \textbf{120} (1998), 955--980.

\bibitem{keraani}
S.~Keraani, \emph{On the blow up phenomenon of the critical nonlinear Schr\"odinger equation}, J. Funct. Anal. \textbf{235} (2006), 171--192.

\bibitem{kwong}
M. K. Kwong, \emph{Uniqueness of positive solutions of $\Delta u-u+u^p=0$ in $\R^n$, }Arch. Rat. Mech. Anal. \textbf{105} (1989),
243-266.

\bibitem{mcleodserrin} K. Mcleod, J. Serrin, \emph{Nonlinear Schr\"odinger equation. Uniqueness of positive solutions
of $\Delta u+f(u)=0$ in $\R^n$}, Arch. Rat. Mech. Anal. \textbf{99} (1987), 115-145.

\bibitem{merle1}
F. Merle, \emph{Determination of blow-up solutions with minimal mass for nonlinear Schr\"odinger equations
with critical power}, Duke Math. J. \textbf{69} (1993), 427--454.

\bibitem{merle-raphael1}
F. Merle, P. Raphael, \emph{The blow-up dynamic and upper bound on the blow-up rate for critical
nonlinear Schr\"odinger equation}, Ann. of Math. \textbf{161} (2005), 157--222.

\bibitem{merle-raphael2}
F. Merle, P. Raphael, \emph{Profiles and quantization of the blow up mass for critical nonlinear
Schr\"odinger equation}, Comm. Math. Phys. \textbf{253} (2005), 675--704.

\bibitem{merle-tsutsumi}
F. Merle, Y. Tsutsumi, \emph{$L^2$-concentration of blow-up
solutions for the nonlinear Schr\"odinger equation with the critical
power nonlinearity}, J. Diff. Eq. \textbf{84} (1990), 205-214.

\bibitem{staf}
G. Staffilani, \emph{On the generalized Korteweg-de Vries-type equations}, Differential Integral Equations \textbf{10} (1997), 777--796.

\bibitem{strichartz}
R.~S.~Strichartz, \emph{Restrictions of Fourier transforms to quadratic surfaces and decay of solutions of wave equations},
Duke Math. J. \textbf{44} (1977), 705--714.

\bibitem{tao:focusing}
T. Tao, \emph{On the asymptotic behavior of large radial data for a focusing non-linear Schr\"odinger
equation,} Dyn. PDE \textbf{1} (2004), 1-48.

\bibitem{taylor}
M. E. Taylor, \emph{Tools for PDE. Pseudodifferential operators, paradifferential operators, and layer potentials},
Mathematical Surveys and Monographs, \textbf{81}, American Mathematical Society, Providence, RI, 2000.

\bibitem{tzirakis} N. Tzirakis, \emph{Mass concentration phenomenon for the quintic nonlinear Schr\"odinger equation in one dimension},
SIAM J. Math. Anal. \textbf{37} (2006), 1923--1946.

\bibitem{Monica:thesis} M. Visan, \emph{The defocusing energy-critical nonlinear Schr\"odinger equation
in dimensions five and higher,} Ph.D. Thesis, UCLA.

\bibitem{MV} M. Visan, \emph{The defocusing energy-critical nonlinear Schr\"odinger equation in higher dimensions,} to appear Duke Math. J.

\bibitem{weinstein:ground state}
M.~I.~Weinstein, \emph{Nonlinear Schr\"odinger equations and sharp interpolation estimates}, Comm. Math. Phys. \textbf{87} (1982/83), 567--576.

\bibitem{weinstein:M} M.~I.~Weinstein, \emph{Modulational stability of ground states of nonlinear Schr\"odinger equations,}
SIAM J. Math. Anal. \textbf{16} (1985), 472--491.

\bibitem{weinstein:L} M.~I.~Weinstein, \emph{Lyapunov stability of ground states of nonlinear Schr\"odinger equations,}
CPAM \textbf{39} (1986), 51--68.

\bibitem{zak} V.E. Zakharov, E.A. Kuznetsov, \emph{Quasi-classical theory for three-dimensional wave collapse},
Sov. Phys. JETP \textbf{64} (1986), 773--380.


\end{thebibliography}
\end{document}